\documentclass[11pt]{article}
\usepackage{latexsym} 
\input{amssym.def} 
\input{amssym} 
\newfont{\fra}{eufm10 scaled 1095} 
\newfont{\Bb}{msbm10 scaled 1095} 
\newfont{\Bbg}{msbm10 scaled 1280} 
\newcommand\CC{{{\Bbb C}}} 
\newcommand\RR{{\mbox{\Bb R}}} 
\newcommand\NN{{\mbox{\Bb N}}} 
 
\newcommand\HH{{\mbox{\Bb H}}}

\newcommand\fg{{\frak{g}}} 
\newcommand\fh{{\frak h}} 
\newcommand\fri{{\frak i}} 
\newcommand\fj{{\frak j}} 
\newcommand\fl{{\frak l}}

\newcommand\fa{{\frak a}} 
 
\newcommand\fd{{\frak d}}

\newcommand\fk{{\frak k}} 
\newcommand\fr{{\frak r}} 
\newcommand\fs{{\frak s}} 
 
\newcommand\fz{{\frak z}} 

\newcommand\cC{{\cal C}} 
\newcommand\cZ{{\cal Z}} 
\newcommand\cH{{\cal H}}

\newcommand\ph{\varphi}

\newcommand{\so}{\mathop{{\frak s \frak o}}}
\newcommand{\fsp}{\mathop{{\frak s \frak p}}}

\newcommand{\Aut}{\mathop{{\rm Aut}}}

\newcommand{\SO}{\mathop{{\it SO}}} 
\newcommand{\Hom}{\mathop{{\rm Hom}}} 
\newcommand{\Sp}{\mathop{{{\it Sp}(1)}}} 
\newcommand{\Phl}{\mathop{\Phi_{\fl}}}
\newcommand{\Pha}{\mathop{\Phi_{\fa}}}
\newcommand{\Phg}{\mathop{\Phi_{\fg}}}
\newcommand{\Id}{\mathop{{\rm Id}}} 
 
\newcommand{\ad}{\mathop{{\rm ad}}} 
\newcommand{\tr}{\mathop{{\rm tr}}}

\newcommand{\Ker}{\mathop{{\rm ker}}}

\renewcommand{\Im}{\mathop{{\rm Im}}}

\newcommand{\Span}{\mathop{{\rm span}}} 
\newcommand{\mod}{\mathop{{\rm mod}}} 
\newcommand{\pro}{{\rm pr}} 
\newcommand{\cZQ}{{\cZ^{2}_{Q}(\fl,\fa)}} 
\newcommand{\cHQ}{{\cH^{2}_{Q}(\fl,\fa)}} 
\newcommand{\cCQ}{{\cC^{1}_{Q}(\fl,\fa)}} 
\newcommand{\cZP}{{\cZ^{2}_{Q}(\fl,\Phl,\fa)}} 
\newcommand{\cHP}{{\cH^{2}_{Q}(\fl,\Phl,\fa)}} 
\newcommand{\cCP}{{\cC^{1}_{Q}(\fl,\Phl,\fa)}} 
\newcommand\ip{{\langle\cdot \,,\cdot \rangle}} 
 
\newcommand\lb{{[\cdot\,,\cdot]}} 
 
\newcommand\dd{\fd_{\alpha,\gamma}(\fl,\Phl,\fa)} 
\newcommand\proof{{\sl Proof. }} 
\newcommand{\qed}{\hspace*{\fill}\hbox{$\Box$}\vspace{2ex}} 
\newtheorem{theo}{Theorem}[section] 
\newtheorem{pr}[theo]{Proposition} 
\newtheorem{de}[theo]{Definition} 
 
\newtheorem{re}[theo]{Remark} 
\newtheorem{co}[theo]{Corollary} 
\newtheorem{lm}[theo]{Lemma} 
\begin{document} 
\title{New examples of indefinite hyper-K\"ahler symmetric spaces}
\author{Ines Kath and Martin Olbrich} 
\maketitle 
\begin{abstract}
\noindent
Following the approach to pseudo-Riemannian symmetric spaces developed in
\cite{KO3} we exhibit
examples of indefinite hyper-K\"ahler symmetric spaces
with non-abelian holonomy. Moreover, we classify indecomposable hyper-K\"ahler symmetric
spaces whose metric has signature $(4,4n)$. Such spaces exist if and only if
$n\in\{0,1,3\}$.
\end{abstract}
\section{Introduction}

The theory of 
special pseudo-Riemannian geometries
has been steadily developing for some years now.
Hyper-K\"ahler geometry is certainly one of the most important of these geometries. Recall that a (pseudo-Riemannian) hyper-K\"ahler manifold is a tuple
$(M,g,I,J)$, where $(M,g)$ is a pseudo-Riemannian manifold, 
and $I,J$ are two parallel anticommuting almost complex structures on
$M$ which preserve the scalar product $g$. In particular, $\dim M$ and
the index of $g$ are
divisible by $4$. It is natural to look first
for symmetric examples of such manifolds, i.e.~hyper-K\"ahler
symmetric spaces, and to try to classify them.
It is well-known that
there are no non-flat Riemannian hyper-K\"ahler symmetric spaces. However, the pseudo-Riemannian situation is different.
Alekseevsky and Cort\'es \cite{AC} give a nice construction of
hyper-K\"ahler symmetric spaces of real dimension $8n$ in terms of
quartic polynomials on a $2n$-dimensional complex
vector space
which satisfy a certain reality condition. All these hyper-K\"ahler symmetric spaces have
neutral signature and an abelian holonomy group of very special
structure. A little more detail on this
construction will be given in Section \ref{SEx}, Remark \ref{pol}.

In \cite{KO3} we developed a systematic approach to the construction and
classification of pseudo-Riemannian symmetric spaces. Here we will
specialize  this approach to the hyper-K\"ahler case.
It is based on the fact that the Lie
algebra of the transvection group of a symmetric space can be obtained from simpler objects by a canonical extension procedure, called quadratic extension.
It is the aim of the present paper to shed some new light to the theory
of hyper-K\"ahler symmetric spaces from this rather different
perspective.

Simply connected hyper-K\"ahler symmetric spaces can be
described by so-called hyper-K\"ahler symmetric triples. Such a triple consists
of the Lie algebra of the transvection group of the symmetric space, an
invariant non-degenerate inner product on this Lie algebra and a certain
$\Sp$-action which is called quaternionic grading, see Section~\ref{SHyp}
for an exact definition. In Section~\ref{SCon} we will describe an extension procedure
which yields a hyper-K\"ahler symmetric triple starting with a Lie algebra
with quaternionic grading $(\fl,\Phl)$, a pseudo-Euclidean vector space $\fa$ with
quaternionic grading $\Pha$ and a pair $(\alpha,\gamma)\in 
\Hom(\bigwedge^{2}\fl,\fa)\oplus\Hom(\bigwedge^{3}\fl,\RR)$ of
$\Sp$-invariant forms satisfying
certain cocycle conditions.
Note that the transvection group of a Hyper-K\"ahler symmetric space
is always nilpotent, see Corollary \ref{nico}.
Therefore the same is true for the Lie algebra $\fl$.

In Section \ref{SEx} we will use this general construction to give concrete examples. We construct
hyper-K\"ahler symmetric triples whose associated symmetric spaces
have a metric of signature $(4n+4,4n+12)$, $n\ge 0$, and a non-abelian
holonomy group. Using the tangent bundle construction we obtain
further examples which also have non-abelian holonomy.
All these examples are indecomposable, i.e.~they cannot be written as the direct
sum of triples of smaller dimension.
Thus, there are much more hyper-K\"ahler symmetric spaces than those
which arise by the construction due to Alekseevsky and
Cort\'es \cite{AC} mentioned at the beginning. In particular, the existence of these examples contradicts the classification result for hyper-K\"ahler symmetric spaces in \cite{AC} and also has consequences for the results in \cite{ABC} and \cite{DJS}. For more information see Remark \ref{pol}.

Furthermore, we show that there is a canonical way to represent each
hyper-K\"ahler symmetric triple as an extension of the kind described above.
Such extensions of $(\fl,\Phl)$ by $\fa$ which
are associated with a hyper-K\"ahler symmetric space in this canonical way are
classified by a cohomology set $\cHP_{\sharp}$ (introduced in Section \ref{SQua}). The isomorphism
classes of underlying hyper-K\"ahler symmetric triples are in 
correspondence with the orbits of the action of the product of
the automorphism groups of $(\fl,\Phl)$ and $\fa$ on $\cHP_{\sharp}$. 
This yields a general classification scheme for hyper-K\"ahler symmetric
triples. Theorem \ref{general} states the corresponding result for
indecomposable triples, i.e. for the objects one really wants to classify.

This general classification scheme can be used to find explicit
classification results (i.e.~lists) if one considers only
hyper-K\"ahler symmetric spaces with a metric of a given small index.
In Section \ref{S4} we will demonstrate this for indecomposable spaces
of index~4. They are exhausted by the flat space $\HH$, a one-parameter family of hyper-K\"ahler symmetric spaces of signature $(4,4)$ obtainable by the Alekseevsky-Cort\'es construction, and one single space of signature $(4,12)$ (the one constructed in Section~\ref{SEx}, Example 1). See Theorem \ref{T4} for the precise statement.

\section{Hyper-K\"ahler symmetric triples} \label{SHyp}

We will say that $(\fl,\Phi_{\fl})$ is a {\it Lie algebra with quaternionic
grading}, if $\fl$ is a Lie algebra and
$\Phi_{\fl}:\Sp\rightarrow \Aut(\fl)$ is an $\Sp$-action
of the following kind. We assume that $\fl=\fl_{+}\oplus\fl_{-}$
such that the representation of $\Sp$ on $\fl_{+}$ is trivial and
the representation of $\Sp$ on $\fl_{-}$ is a multiple of the standard
representation. In particular, $\fl_{-}$ is a left $\HH$-module.
If in addition $[\fl_{-},\fl_{-}]=\fl_{+}$ we call the grading {\it proper}.

Similarly, we will say that $(\fa,\ip_\fa,\Pha)$ (or $\fa$ in abbreviated 
notation) is a {\it vector space with 
orthogonal quaternionic grading}, if $(\fa,\Pha)$ is an abelian Lie 
algebra with quaternionic grading and the image of $\Pha$ is in 
$O(\fa,\ip_{\fa})$. 

\begin{pr}\label{nilpferd} 
If $(\fl,\Phi_{\fl})$ is a Lie algebra with proper quaternionic grading, 
then $\fl$ is nilpotent. 
\end{pr} 
\proof 
We consider the semi-simple Lie algebra $\fs=\fl/\fr$, where $\fr$ 
is the solvable radical of $\fl$. The grading $\Phi_\fl$ induces a proper 
quaternionic grading $\Phi_\fs$. 
Being a connected subgroup of the automorphism group of $\fs$, the image 
of $\Phi_\fs$ consists of {\em inner} automorphisms 
which respect the decomposition $\fs=\fs_{+}\oplus\fs_{-}$. Therefore 
its Lie algebra can be identified with a subalgebra $\fk\subset \fs_+$. 
It follows that the adjoint representation of $\fk$ on $\fs_+$, and hence on 
$\fk$, is trivial. On the other hand, $\fk\cong \frak{sp}(1)$ unless 
$\fs_{-}=\{0\}$. 
We conclude that $\fs=\{0\}$, i.e., $\fl$ is solvable.

We finish the proof by showing that solvable Lie algebras 
$\fl=\fl_+\oplus\fl_{-}$, $\fl_+=[\fl_-,\fl_-]$, admitting an automorphism $F$ 
such that 
$F^2|_{\fl_-}=-\Id$ and $F|_{\fl_+}=\Id$ are nilpotent. 

We look at the decreasing chain of ideals of $\fl$ 
$$ \fl:=R_0(\fl)\supset R_1(\fl)\supset R_2(\fl)\dots \ ,$$ 
where $R_{k+1}(\fl)$ is the minimal $\fl$-ideal in $R_k(\fl)$ such that 
the induced action of $\fl$ on $R_k(\fl)/R_{k+1}(\fl)$ is semi-simple. 
Note that these ideals 
are invariant under all automorphisms of $\fl$. 
There exists a number $m$ such that $R_m(\fl)=\{0\}$. 

We consider the automorphism $F=\Phi_\fl(i)$. 
We look at the semi-simple representation of $\fl$ 
on the complexification $W_k$ of $R_k(\fl)/R_{k+1}(\fl)$. 
Since $\fl'$ is a nilpotent ideal of $\fl$, the Lie algebra 
$\fl_+=[\fl_-,\fl_-]\subset \fl'$ acts trivially on $W_k$. 
It follows that $W_k$ is the direct sum of weight spaces $E_\lambda$ 
with $\lambda\in (\fl_-^*)_\CC$. Note that $F$ acts naturally on $(\fl_-^*)_\CC$ 
and on $W_k$ with the property 
$ F(E_\lambda)=E_{F(\lambda)}$. Assume that $\lambda\ne 0$.
Then the elements
$F^n(\lambda)$, $n=0,1,2,3$, are pairwise different. Therefore
the sum of the weight
spaces $E_{F^n(\lambda)}$, $n=0,1,2,3$, is direct. Take $v\in E_\lambda$.
Then $F^{n}(v)\in E_{F^n(\lambda)}$, 
and 
$$v^-:=v-F(v)+F^{2}(v)-F^{3}(v)$$ 
satisfies $F(v^-)=-v^-$. 
However, the only possible 
eigenvalues of $F$ on $W_k$ are $1$, $i$ and $-i$. We conclude that 
$v^-=0$, hence $v=0$. It follows that $W_k=E_0$, i.e. $\fl$ acts trivially on $W_k$. We conclude that $\fl$ is nilpotent.
\qed

\begin{de} 
A hyper-K\"ahler symmetric triple is a triple $(\fg,\Phi_\fg,\ip_\fg)$,
where $(\fg,\Phi_\fg)$ is a Lie algebra with proper quaternionic grading 
and $\ip_\fg$ is an $\fl-$ and $\Phi_\fl$-invariant non-degenerate symmetric 
bilinear 
symmetric bilinear form. 
\end{de} 

There is an obvious notion of isomorphism between hyper-K\"ahler symmetric 
triples. Note that $\theta_\fg:=\Phi_\fg(-1)$ is an isometric involution of 
$\fg$ and $(\fg,\theta_\fg,\ip_\fg)$ is a symmetric triple in the sense of 
\cite{KO3}. 

\begin{pr} 
The Lie algebra of the transvection group of a hyper-K\"ahler symmetric space carries the structure of a hyper-K\"ahler symmetric triple in a canonical way.
There is a one-to-one correspondence between isometry classes of simply connected
hyper-K\"ahler symmetric spaces and isomorphism classes of hyper-K\"ahler 
symmetric triples. 
\end{pr}
\proof
The proposition is a slight variant of the well-known correspondence between
pseudo-Riemannian symmetric spaces and symmetric triples (see \cite{CP80}, 
compare \cite{KO3}, Section 2). If $(\fg,\Phi_\fg,\ip_\fg)$ is a hyper-K\"ahler 
symmetric triple and $(M,g)$ is a symmetric space with associated symmetric 
triple 
$(\fg,\theta_\fg,\ip_\fg)$, then $I=\Phi_\fg(i)|_{\fg_-}$ and 
$J=\Phi_\fg(j)|_{\fg_-}$ 
are $\fl_+$-invariant anticommuting complex structures on $\fg_-\cong T_oM$ 
respecting the metric and therefore 
induce a hyper-K\"ahler structure on $M$. 

The only non-obvious point is the opposite direction. Let $(M,g,I,J)$ be a hyper- 
K\"ahler symmetric space and let $(\fg,\theta_\fg,\ip_\fg)$ be the associated
symmetric triple. Then $I,J,K:=IJ$ span a Lie algebra $\fk\cong \frak{sp}(1)$ 
which acts orthogonally on $\fg_-\cong T_oM$. This action commutes with the one 
of $\fg_+$. We extend the $\fk$-action to $\fg$ by the trivial action on $\fg_+$. 
For $X,Y\in\fg_-$, $Z\in\fg_+$, and $Q\in\fk$ we compute 
\begin{eqnarray*} 
\langle [QX,Y]+[X,QY], Z\rangle &=& 
\langle QX,[Y,Z]\rangle +\langle X,[QY,Z]\rangle \\ 
&=& 
\langle QX,[Y,Z]\rangle +\langle X,Q[Y,Z]\rangle=0\ . 
\end{eqnarray*} 
It follows that $\fk$ acts by derivations on $\fg$. Integrating the resulting 
homomorphism of $\frak{sp}(1)$ into the derivations of $\fg$ we obtain the 
desired 
homomorphism $\Phi_\fg: \Sp\rightarrow \Aut(\fg)$. 
\qed 

Taking Proposition \ref{nilpferd} into account we obtain

\begin{co}\label{nico}
The transvection group of a hyper-K\"ahler symmetric space is nilpotent.
\end{co}

A hyper-K\"ahler symmetric triple is called decomposable, if it is the orthogonal
direct sum of two non-zero hyper-K\"ahler symmetric triples, and indecomposable
otherwise. The simply connected hyper-K\"ahler symmetric spaces
which correspond to indecomposable hyper-K\"ahler triples are precisely those
which are indecomposable in the differential geometric sense. This follows
from the de Rham--Wu decomposition theorem.

The underlying metric Lie algebra of a hyper-K\"ahler symmetric triple 
$(\fg,\Phi_\fg,\ip_\fg)$ is the tuple $(\fg,\ip_\fg)$. There are
analogous notions of indecomposability for symmetric triples and for metric Lie
algebras.

\begin{lm}\label{A0}
Let $(\fg,\Phi_\fg,\ip_\fg)$ be a hyper-K\"ahler symmetric triple
with $\fg$ non-abelian.
Then the following conditions are equivalent:
\begin{enumerate}
\item[(i)]$(\fg,\Phi_\fg,\ip_\fg)$ is indecomposable.
\item[(ii)] The symmetric triple $(\fg,\theta_\fg,\ip_\fg)$ is indecomposable.
\item[(iii)] The metric Lie algebra $(\fg,\ip_\fg)$ is indecomposable.
\end{enumerate}
\end{lm}
\proof
The implications $(iii)\Rightarrow (ii)\Rightarrow (i)$ are obvious.
We have to prove $(i)\Rightarrow (iii)$.

Assume $(i)$. We first observe that
the center $\fz(\fg)$ is isotropic, i.e. $\fz(\fg)\subset \fz(\fg)^\perp=\fg'$.
Indeed, any $\Sp$-invariant complement of $\fz(\fg)\cap \fz(\fg)^\perp$ in $\fz(\fg)$
is a non-degenerate abelian ideal of $\fg$. It has to be zero by assumption.
Let
\begin{equation}\label{deco}
\fg=\bigoplus_{k=1}^n \fri_k
\end{equation}
be an orthogonal decomposition of $(\fg,\ip_\fg)$
into non-trivial indecomposable ideals.
Let $\fj\subset \fg$ be a further non-degenerate indecomposable ideal. According
to \cite{A78}, Theorem~3, there exists $l\in\{1,\dots,n\}$ such that
$\fj^\prime={\fri_l}'$.
In particular, the set of
derivatives of non-degenerate indecomposable ideals is finite. Therefore the
action of the connected group $\Sp$ on this set is trivial. It follows
that the ideals
${\fri_1}'$ and $(\bigoplus_{k=2}^n \fri_k')^\perp=\fri_1+\fz(\fg)$ are $\Sp$-invariant. The above observation
concerning the center yields $\fz(\fri_1)={\fri_1}'\cap\fz(\fg)$.
It follows that this ideal is $\Sp$-invariant, too. Therefore $\ip_\fg$ induces an $\Sp$-equivariant surjective map
$$ \psi:\fri_1+\fz(\fg)\longrightarrow \fz(\fri_1)^*\ .$$
Note that ${\fri_1}'\subset \Ker(\psi)$. Let $s:\fz(\fri_1)^*\rightarrow\fri_1+\fz(\fg)$ be an $\Sp$-equivariant
section of $\psi$. Then $\widetilde\fri_1:={\fri_1}'\oplus s(\fz(\fri_1)^*)$ is a non-degenerate and $\Sp$-invariant ideal of $\fg$.
Moreover, the  projection $\widetilde\fri_1\rightarrow\fri_1$ with respect to
(\ref{deco}) is an isomorphism.
Indecomposability of the triple $(\fg,\Phi_\fg,\ip_\fg)$ now implies $\fg=\widetilde\fri_1$. Thus $\fg\cong \fri_1$ is indecomposable
as a metric Lie algebra.
\qed

\section{Quaternionic quadratic cohomology}\label{SQua}

Let us first recall the notion of quadratic cohomology introduced in \cite{KO2}.
Let $\fl$ be a finite-dimensional Lie algebra. An {\it orthogonal $\fl$-module} 
is a
tuple $(\rho,\fa,\ip_{\fa})$ (also $\fa$ in abbreviated notation) such that 
$\rho$ is a representation of
the Lie algebra 
$\fl$ on the finite-dimensional real vector space $\fa$ and $\ip_\fa$ is
a scalar product on $\fa$ such that $\rho(L)\in \so(\fa,\ip_\fa)$
for all $L\in\fl$.

For $\fl$ and (any $\fl$-module) $\fa$ we have the
standard cochain complex $(C^*(\fl,\fa),d)$, where $C^p(\fl,\fa)=
\Hom(\bigwedge ^p\fl,\fa)$
and we have the corresponding cocycle groups $Z^p(\fl,\fa)$ and cohomology groups 
$H^p(\fl,\fa)$.
If $\fa$ is the one-dimensional trivial
representation, then we denote this cochain complex also by $C^*(\fl)$.

We have a product
$$\langle\cdot \wedge \cdot\rangle:\, C^p(\fl,\fa)\times C^q(\fl,\fa)\longrightarrow C^{p+q}(\fl)$$
defined by the composition 
$$C^p(\fl,\fa)\times C^q(\fl,\fa)\stackrel{\wedge}{\longrightarrow} 
C^{p+q}(\fl,\fa\otimes \fa) \stackrel{\ip_\fa}{\longrightarrow} C^p(\fl).$$ 

The group of quadratic $1$-cochains is the group
$${\cal C} ^{1}_Q(\fl,\fa)=C^{1}(\fl,\fa)\oplus C^{2}(\fl) $$
with group operation defined by 
$$ (\tau_1,\sigma_1)*(\tau_2,\sigma_2)=(\tau_1+\tau_2, \sigma_1 
+\sigma_2 +\textstyle{\frac12} \langle \tau_1\wedge \tau_2\rangle)\,.$$ 
We consider now the set  
$${\cal Z} ^{2}_Q(\fl,\fa)=\{(\alpha,\gamma) \in C^{2}(\fl,\fa)\oplus 
C^{3}(\fl) \mid d\alpha=0,\ 
d\gamma=\textstyle{\frac12}\langle\alpha \wedge\alpha\rangle\} $$ 
whose elements are called quadratic $2$-cocycles. The group  
${\cal C} ^{1}_Q(\fl,\fa)$ acts on ${\cal Z} ^{2}_Q(\fl,\fa)$ by 
$$(\alpha,\gamma)(\tau,\sigma)=\Big(\,\alpha +d\tau,\gamma +d\sigma 
+\langle(\alpha +\textstyle{\frac12} d\tau)\wedge\tau\rangle\,\Big).$$
Ordinary quadratic cohomology is then the orbit space of this action:
$$
\cHQ=\cZQ/\cCQ\ .$$
Now let $(\fl,\Phl)$ be a Lie algebra with quaternionic grading and
let $(\fa,\Pha)$ be a vector space with orthogonal quaternionic
grading. We consider $\fa$ as a trivial $\fl$-module. Then $\Phi_{\fa}$
and $\Phl$ define $\Sp$-actions on $\cZQ$ and $\cCQ$. More precisely, for
$q\in\Sp$ the pair of morphisms $(\Phl(q),\Pha(q)^{-1})$ induces pullback
maps on $C^{2}(\fl,\fa)\oplus C^{3}(\fl)$ and on $\cCQ=C^{1}(\fl,\fa)\oplus
C^{2}(\fl)$, which leave invariant $\cZQ$ and are compatible with the $\cCQ$-action on $\cZQ$. We consider the
sets of invariants
$$\cZP:={\cZQ}^{\Sp} \ \mbox{ and }\ \cCP:={\cCQ}^{\Sp} \ .$$
The group $\cCP$ acts on $\cZP$ and we can define the (second) quaternionic quadratic
cohomology by 
$$\cHP:=\cZP/\cCP\ .$$

For a Lie algebra $\fl$ we denote by $\fl^1=\fl,\dots, 
\fl^k=[\fl,\fl^{k-1}],\dots$ the lower central series. 
\begin{de}\label{adm}
 Let $(\fl,\Phl)$ be a Lie algebra with a proper quaternionic grading and let 
  $\fa$ be a vector space with 
 orthogonal quaternionic grading. Let 
 $m$ be such 
 that $\fl^{m+2}=0$. Let $(\alpha,\gamma)\in\cZP$. Then the cohomology class 
 $[\alpha,\gamma]\in\cHP$ is called admissible if and only if 
 the following 
 conditions $(T)$, $(A_k)$ and $(B_k)$ hold for all $0\le k\le m$. 
\begin{enumerate} 
\item[$(T)$]$\quad \fa_+=\alpha(\Ker \lb_{\fl_-})$. 
\item[$(A_k)$] 
Let $L_0\in \fz(\fl)\cap \fl^{k+1}$ be such that there exist 
elements 
$A_0\in \fa$ and $Z_0\in (\fl^{k+1})^*$ satisfying 
\begin{enumerate} 
\item[(i)] $\alpha(L,L_0)=0 $, 
\item[(ii)] $\gamma(L,L_0,\cdot)=-\langle A_0,\alpha(L,\cdot)\rangle_\fa 
+\langle 
Z_0, [L,\cdot]_\fl\rangle$ as an element of $(\fl^{k+1})^*$, 
\end{enumerate} 
for all $L\in\fl$, then $L_0=0$. 
\item[$(B_k)$] The subspace $\alpha(\ker \lb_{\fl\otimes \fl^{k+1}})\subset 
\fa$ is non-degenerate, where $\ker \lb_{\fl\otimes \fl^{k+1}}$ is the kernel of 
the map $\lb:\fl\otimes \fl^{k+1}\rightarrow \fl$. 
\end{enumerate} 
We denote the set of all admissible cohomology classes by $\cHP_{\sharp}$. 
\end{de} 
The admissibility conditions are specializations of the ones in \cite{KO3}, 
Definition 5.2, 
to the case that $\fl$ is nilpotent and that the representation of 
$\fl$ on $\fa$ is trivial. As in \cite{KO3}, they do not depend
on the choice of the cocycle representing the cohomology class $[\alpha,\gamma]$.

Now let $\fl_i$, $\fa_i$, $i=1,2$, be Lie algebras (vector spaces, resp.)
with (orthogonal) quaternionic grading. We form $\fl=\fl_1\oplus \fl_2$,
$\fa=\fa_1\oplus\fa_2$ and consider the corresponding projections 
$\pro_i:\fl\rightarrow \fl_i$. If 
$(\alpha_i,\gamma_i)\in\cZ^{2}_{Q}(\fl_i,\Phi_{l_i},\fa_i)$, then 
$$ (\pro_1^*\alpha_1\oplus 
\pro_2^*\alpha_2,\pro_1^*\gamma_1+\pro_2^*\gamma_2)\in\cZP\ .$$ 
This operation induces a map 
$$\oplus:\cH^{2}_{Q}(\fl_1,\Phi_{l_1},\fa_1) 
\oplus\cH^{2}_{Q}(\fl_2,\Phi_{l_2},\fa_2)
\longrightarrow \cHP\ .$$ 

\begin{de} 
A cohomology class $\phi\in\cHP$ is called decomposable, if there are 
a decomposition $\fl=\fl_1\oplus \fl_2$ into $\Phi_\fl$-invariant ideals, a 
$\Phi_\fa$-invariant orthogonal decomposition $\fa=\fa_1\oplus\fa_2$ 
(at least one of these decompositions should be non-trivial), and 
cohomology classes $\phi_i\in\cH^{2}_{Q}(\fl_i,\Phi_{l_i},\fa_i)$ such that 
$\phi=\phi_1\oplus \phi_2$. 
Otherwise we call $\phi$ indecomposable. 

We denote the set of all indecomposable elements in $\cHP_{\sharp}$ by 
$\cHP_{0}$. 
\end{de} 

\section{A construction method}\label{SCon}
Now we will present a construction method which yields hyper-K\"ahler symmetric
triples starting with a Lie algebra with proper quaternionic grading 
$(\fl,\Phl)$, a vector space $\fa$ with orthogonal quaternionic grading and an 
admissible 
cocycle in $\cZP$ (i.e. a cocycle which represents an admissible cohomology 
class). This construction is a special case of the construction method for 
symmetric triples presented in \cite{KO3}, Section 4.2. In Section \ref{SClass} 
we will see that each hyper-K\"ahler symmetric triple arises by this 
construction in a canonical way. 

Let $(\fl,\Phl)$ be a Lie algebra with proper quaternionic grading and let 
$(\fa,\ip_\fa,\Pha)$ be a vector space with orthogonal quaternionic grading. We 
consider the vector space $$\fd:=\fl^*\oplus\fa\oplus\fl.$$ 
Now we choose 
an admissible cocycle 
$(\alpha,\gamma) \in \cZP$ and define a Lie bracket 
$\lb:\fd\times\fd\rightarrow \fd$ by $[\fl^*\oplus\fa,\fl^*\oplus\fa] 
=0$ and 
\begin{eqnarray} \label{lb1}
\ [L_1,L_2] &=& \gamma(L_1,L_2,\cdot) +\alpha(L_1,L_2)+[L_1,L_2]_\fl\\ 
\ [L,A] &=& - \langle A,\alpha(L,\cdot)\rangle\\ 
\ [L,Z]& = & \ad{}^*(L)(Z) \label{lb3}
\end{eqnarray} 
for $Z\in \fl^*$, $A\in \fa$ and 
$L,\,L_1,\,L_2\in \fl$. 
Moreover we define 
an inner product $\ip$ and a quaternionic grading $\Phi$ on $\fd$ by 
\begin{eqnarray*} 
\langle Z_1+A_1+L_1,Z_2+A_2+L_2\rangle&:=& \langle A_1,A_2\rangle_\fa 
+Z_1(L_2) +Z_2(L_1) \\ 
\Phi(Z+A+L)&:=& \Phi_{\fl}^*(Z)+\Pha(A)+\Phl(L) 
\end{eqnarray*} 
for $Z,\,Z_1,\,Z_2\in \fl^*$, $A,\,A_1,\,A_2\in \fa$ and 
$L,\,L_1,\,L_2\in \fl$. 

\begin{pr} \label{P41}
 The tuple $\dd:=(\fd,\Phi,\ip)$ is a hyper-K\"ahler symmetric triple. 
 It is indecomposable if and only if $[\alpha,\gamma]\in\cHP_\sharp$ 
 is indecomposable. 
\end{pr} 
\proof 
This follows from the results of \cite{KO3}, especially Prop.~4.1, Lemma~5.2, and Prop.~6.2.
\qed
\begin{re}
{\rm
If $(\alpha,\gamma)\in \cZP$ is arbitrary, then $\dd$ is still a Lie algebra with quaternionic grading and an invariant inner product. One needs some conditions
on the cocycle in order to ensure properness. For this purpose, conditions
much weaker than admissibility would suffice (e.g. $(T)$ together with $(A_0)$).
However, we can detect indecomposability of $\dd$ for admissible cocycles $(\alpha,\gamma)$ only. The main reason for considering the admissibility conditions is
that any hyper-K\"ahler symmetric triple is isomorphic to $\dd$
for an essentially unique tuple $(\fl,\Phl,\fa,[\alpha,\gamma]\in\cHP_\sharp)$, see Theorem \ref{general}
 below.}
\end{re}

\begin{re}
{\rm
The signature $(p,q)$ of a hyper-K\"ahler symmetric triple 
$(\fg,\Phi_\fg,\ip_\fg)$ is defined to be the signature of the restriction of 
$\ip_\fg$ to $\fg_-$. Here $p$ is the dimension of a maximal negativ definite 
subspace of $\fg_-$. Sometimes $p$ is called index of the triple. The signature 
(the index) of a hyper-K\"ahler symmetric triple equals the signature (the 
index) of the metric of any pseudo-Riemannian (hyper-K\"ahler)
symmetric space
which is associated with the triple. The signature $(p,q)$ of the above 
constructed hyper-K\"ahler symmetric triple $\dd$ is determined by the 
signature $(p_\fa,q_\fa)$ of the restriction of $\ip_\fa$ to $\fa_-$ and the 
dimension of $\fl_-$: 
\begin{equation}\label{signature} 
(p,q)=(p_\fa+\dim \fl_-,q_\fa+\dim \fl_-) 
\end{equation} 
} 
\end{re}

%
\section{Examples of hyper-K\"ahler symmetric spaces with 
non-abelian holonomy} \label{SEx} 
Now we will use the method described above to construct various hyper-K\"ahler 
symmetric triples whose associated symmetric spaces have non-abelian holonomy. 
In particular, these spaces are missing in the classification result for 
hyper-K\"ahler symmetric triples claimed by Alekseevsky and Cort\'es \cite{AC}. 

\subsection*{Example 1} 
First let us define the structure of a Lie algebra with proper quaternionic 
grading on the vector space 
$\fl_{0}=\HH\oplus \Im \HH$. For all $q\in\HH$ we 
denote $(q,0)$ only by $q$ and $(0,i)$, $(0,j)$, $(0,k)$ by $I$, $J$, 
and $K$ respectively. On $\fl_{0}$ we define a Lie bracket $\lb$ 
by 
$$I,J,K\in \fz(\fl_{0}),\qquad 
[q_{1},q_{2}]=(0,\Im \bar q_{1} q_{2})\in \HH\oplus \Im \HH.$$ 
Moreover, on $\fl_{0}=\HH\oplus \Im \HH$ we consider the $\Sp$-action 
$\Phi_{\fl_0}$ which is 
the left multiplication on the first summand and which is trivial on 
the second summand. 

Now we will describe a suitable vector space $\fa_0$ with orthogonal
quaternionic grading such that $\cH_Q^2(\fl_0,\Phi_{\fl_0},\fa_{0})_0$ is not 
empty. 
Let $e_{1}, e_{2}, e_{3}$ be the standard basis of $\RR^{3}$ and 
define $A_{1}=e_{1}- e_{2}$, $A_{2}=e_{2} - e_{3}$ and $A_{3}=e_{3} - 
e_{1}$. We consider the 2-dimensional vector space 
$\fa_{\Bbb R}=\Span\{A_{1},A_{2},A_{3}\}.$ Let $\ip$ be the restriction 
of the (positive definite) standard scalar product on $\RR^{3}$ to 
$\fa_{\Bbb R}$ and let $\ip_{\Bbb H}$ be the standard scalar product 
on $\HH$, i.e. $\langle p,q\rangle_{\Bbb H} = \bar pq$. Then 
$\ip_{0}:=\ip_{\Bbb H}\otimes \ip$ is a scalar product on 
$\fa_{0}:=\HH\otimes_{\Bbb R}\fa_{\Bbb R}$.
Furthermore, if we let $\Phi_{\fa_{0}}$ be the left multiplication on 
$\fa_{0}$ we obtain a vector space with 
orthogonal quaternionic grading $(\fa_{0},\ip_0,\Phi_{\fa_{0}})$. 

We define $\alpha_0\in C^{2}(\fl_{0},\fa_{0})$ by 
\begin{equation}\label{alpha0} 
\alpha_0(q,I)=qiA_{1}, \quad \alpha_0(q,J)=qjA_{2}, 
\quad\alpha_0(q,K)=qkA_{3} 
\end{equation} 
for all $q\in \HH$ and by 
$\alpha_0((\fl_{0})_{+},(\fl_{0})_{+})=\alpha_0((\fl_{0})_{-},(\fl_{0})_{-})=0$. 
Furthermore, we define $\gamma_0\in C^{3}(\fl_{0})$ by 
$\gamma_{0}(I,J,K)=2$ and $\gamma_{0}((\fl_{0})_{-},\fl_{0},\fl_{0})=0$.

\begin{lm}\label{lemma} 
 We have $(\alpha_0,\gamma_{0})\in \cZ^{2}_{Q}(\fl_{0},\Phi_{\fl_{0}},\fa_{0})$.
 Moreover, $(\alpha_0,\gamma_{0})$ is admissible and indecomposable.
\end{lm}
\proof Obviously, $\alpha_0$ and $\gamma_0$ are $\Sp$-invariant. Let us compute
\begin{eqnarray*} 
 d\alpha_0(1,i,j) 
 &=&-\alpha_0([1,i],j)+\alpha_0([1,j],i)-\alpha_0([i,j],1)\\ 
 &=&-\alpha_0(I,j)+\alpha_0(J,i)+\alpha_0(K,1)\\ 
 &=& jiA_{1}-ijA_{2}-kA_{3}\ =\ 0. 
\end{eqnarray*} 
Because of the $\Sp$-invariance of $\alpha_{0}$ this implies 
$$ d\alpha_0(1,i,k) =d\alpha_0(1,j,k)=d\alpha_0(i,j,k)=0.$$ 
Since all other components of $d\alpha_0$ vanish obviously, we obtain 
$d\alpha_0=0$. 

Now we will check the condition $d\gamma_{0}=\frac 
12\langle\alpha_0\wedge\alpha_0\rangle$. Obviously we have 
$$d\gamma_{0}(\fl_0,(\fl_0)_{+},(\fl_0)_{+},(\fl_0)_{+})= 
\langle\alpha_0\wedge\alpha_0\rangle(\fl_0,(\fl_0)_{+},(\fl_0)_{+},(\fl_0)_{+})=0
$$ 
and 
$$d\gamma_{0}(\fl_0,(\fl_0)_{-},(\fl_0)_{-},(\fl_0)_{-})= 
\langle\alpha_0\wedge\alpha_0\rangle 
(\fl_0,(\fl_0)_{-},(\fl_0)_{-},(\fl_0)_{-})=0.$$ By 
$\Sp$-invariance of $(\alpha_0,\gamma_{0})$ it remains to prove 
$$d\gamma_{0}(1,q,P,Q)=\textstyle{\frac 
12}\langle\alpha_0\wedge\alpha_0\rangle(1,q,P,Q)$$ 
for all imaginary $q\in\HH$ and all $P,Q\in\{I,J,K\}$. 
We will show this for $q=i$. The remaining equations can be proved 
similarly. We have 
\begin{eqnarray*} 
 \textstyle{\frac 12}\langle\alpha_0\wedge\alpha_0\rangle(1,i,I,J)&=&        
  \langle 
 \alpha_0(I,1),\alpha_0(i,J)\rangle_{0}+\langle 
 \alpha_0(i,I),\alpha_0(1,J)\rangle_{0}\\ 
 &=&\langle -iA_{1},kA_{2} \rangle_{0}+\langle -A_{1},jA_{2} 
 \rangle_{0}\\ 
 &=&0\ =\ -\gamma_{0}(I,I,J)\ =\ d\gamma_{0}(1,i,I,J)\\[1ex] 
 \textstyle{\frac 12}\langle\alpha_0\wedge\alpha_0\rangle(1,i,I,K)&=& 
 \langle \alpha_0(I,1),\alpha_0(i,K)\rangle_{0}+\langle 
 \alpha_0(i,I),\alpha_0(1,K)\rangle_{0}\\ 
 &=&\langle -iA_{1},-jA_{3} \rangle_{0}+\langle -A_{1},kA_{3} 
 \rangle_{0}\\ 
 &=&0 \ =\ -\gamma_{0}(I,I,K)\ =\ d\gamma_{0}(1,i,I,K)\\[1ex] 
 \textstyle{\frac 12}\langle\alpha_0\wedge\alpha_0\rangle(1,i,J,K)&=&\langle 
 \alpha_0(J,1),\alpha_0(i,K)\rangle_{0}+\langle 
 \alpha_0(i,J),\alpha_0(1,K)\rangle_{0}\\ 
 &=&\langle -jA_{2},-jA_{3} \rangle_{0}+\langle kA_{2},kA_{3} 
 \rangle_{0}\\ 
 &=&- 2\ 
 =\ -\gamma_{0}(I,J,K)\ =\ d\gamma_{0}(1,i,J,K). 
\end{eqnarray*} 
It is easy to see that $(\alpha_0,\gamma_0)$ is
admissible. Indeed, Condition $(T)$ is satisfied because of 
$(\fa_0)_{+}=0$ 
and $(B_0)$ and $(B_1)$ hold since $\ip_0$ is definite. As for Conditions 
$(A_0)$ and $(A_1)$ we observe that $\fz(\fl_0)\cap\fl_0 
=\fz(\fl_0)\cap\fl_0'=\Span\{I,J,K\}$ and that for all $Q\in\Span\{I,J,K\}$ 
there is an element $L\in\fl_0$ such that $\alpha_{0}(L,Q)\not=0$. Hence, $(A_0)$ 
and $(A_1)$ are also satisfied, thus $(\alpha_{0},\gamma_{0})$ is admissible. 
Obviously, 
$(\alpha_{0},\gamma_{0})$ is also indecomposable since $\fl_0$ is indecomposable 
and 
$\alpha_{0}(\fl_0,\fl_0)=\fa_{0}$. 
\qed 
\begin{co} 
The triple $\fd_{\alpha_0,\gamma_0}(\fl_0,\Phi_{\fl_0},\fa_{0})$ is an
indecomposable
hyper-K\"ahler symmetric triple of signature $(4,12)$. The holonomy group
of a  symmetric space associated with this triple is non-abelian. 
\end{co} 
\proof 
The first statement follows from Proposition \ref{P41}, Lemma 
\ref{lemma} and Equation (\ref{signature}). It remains to prove the 
assertion on the holonomy group. The Lie algebra of this group is 
isomorphic to $\fd_{+}=(\fl_{0}^{*})_{+}\oplus 
(\fa_{0})_{+}\oplus(\fl_{0})_{+}$ with Lie bracket defined by  
(\ref{lb1}) -- (\ref{lb3}). Since $\gamma(I,J,K)\not=0$ this Lie 
algebra is not abelian.
\qed

\begin{re}\label{pol}
{\rm
Let us briefly recall a general construction of hyper-K\"ahler
symmetric triples from \cite{AC}.
Let $(E,\omega)$ be a complex symplectic vector space. Any $S\in S^4E$
defines a complex linear subspace $\fh_S\subset\fsp(E,\omega)\cong S^2 E$ by
$$\fh_S=\Span\{S_{v,w}\in S^2E\:|\: v,w\in E\}\ ,$$
where $S_{v,w}$ is the contraction
of $S$ with $v$ and $w$ defined by the symplectic form $\omega$. 
If
\begin{equation}\label{crux}
S\in (S^4E)^{\fh_S}\ ,
\end{equation}
then $\fh_S\subset\fsp(E,\omega)$ is a Lie subalgebra and, moreover, there is a natural Lie bracket on $\fh_S\oplus (\HH\otimes_\CC E)$.
If there exists a Lagrangian subspace $E_+\subset E$ such that $S\in S^4E_+\subset S^4 E$, then  $S$ is a
solution of
(\ref{crux}). Let us call solutions of this kind tame. If $S$ is tame, then
the Lie algebra $\fh_S$ is abelian.

Let $J$ be a quaternionic structure on $E$ such that $J^*\omega=\bar\omega$.
Then $J$ induces real structures, all denoted by $\tau$, on $S^4E$, $S^2E\cong\fsp(E,\omega)$, and on
$\HH\otimes_\CC E$. If $S\in (S^4 E)^\tau$ satisfies (\ref{crux}), then
the real Lie algebra
$$ \fg_S:=(\fh_S)^\tau\oplus (\HH\otimes_\CC E)^\tau $$
carries a canonical structure of a hyper-K\"ahler symmetric triple.
Moreover, all hyper-K\"ahler symmetric triples arise in this way.
The Alekseevsky-Cort\'es construction mentioned in the introduction
produces a hyper-K\"ahler symmetric triple $\fg_S$ starting with a
tame $\tau$-invariant solution $S$.

Now it is not difficult to show that
$\fd_{\alpha_0,\gamma_0}(\fl_0,\Phi_{\fl_0},\fa_{0})\cong\fg_S$
for
$(E,\omega,J,S)$ as follows: $(E,\omega)$ is $8$-dimensional with standard
symplectic basis $p_1,\dots,p_4,q_1,\dots,q_4$, the quaternionic structure $J$ is characterized by $J(p_1)=p_2$, $J(p_i)=q_i$ for $i=3,4$, and
$$ S=p_1^3q_3+\sqrt{3} p_1^2p_2p_4-\sqrt{3} p_1 p_2^2q_4-p_2^3p_3\ .$$
Let us verify directly that
$S$ satisfies (\ref{crux}).
Using the $\tau$-invariance of $S$ we find that
$\fh_S=V+\tau(V)$, where
\begin{eqnarray*}
V&=&\Span\{S_{q_1,q_1},\, S_{q_1,q_2},\, S_{q_1,q_4},\, S_{q_1,p_3},\, S_{q_1,p_4}\}\\
&=&\Span\{\sqrt{3}p_1q_3+p_2p_4,\, p_1p_4-p_2q_4,\, p_1p_2,\, p_1^2,\, p_2^2\}
\ .
\end{eqnarray*}
It suffices to check that $P(S)=0$ for $P$ running through the above set
of basis elements of $V$. Since $S_{p_1}=S_{p_2}=0$ we immediately see that
$P(S)=0$ if $P$ is one of the last three basis elements. It remains to compute
\begin{eqnarray*}
(\sqrt{3}p_1q_3+p_2p_4)(S)&=&
2(\sqrt{3}p_1 S_{q_3}+p_2 S_{p_4})
\:=\:{\textstyle\frac{1}{2}}
(\sqrt{3}p_1p_2^3+p_2(-\sqrt{3}p_1p_2^2))\:=\:0\ ,\\
(p_1p_4-p_2q_4)(S)&=&
2(p_1 S_{p_4}-p_2 S_{q_4})\:=\:{\textstyle\frac{1}{2}}
(p_1(-\sqrt{3} p_1 p_2^2)-p_2(-\sqrt{3} p_1^2p_2))\:=\:0\ .
\end{eqnarray*}
Equation (\ref{crux}) follows. On the other hand, $S$ is not tame since $\{v\in E\:|\:S_v=0\}=\Span\{p_1,p_2\}$ is only two-dimensional. This shows that
Theorem~3 in \cite{AC} claiming that all solutions of (\ref{crux}) are tame
is not true. In particular, those results of \cite{AC},\cite{ABC},\cite{DJS}
which are based on this theorem have to be reconsidered.
}
\end{re}

\subsection*{Example 2}

We look at $\HH^2$ as an abelian Lie algebra equipped with the
quaternionic grading
given by the left $\HH$-module structure. We equip $\Im \HH\oplus \Im \HH$
with the trivial $\Sp$-action. Let $A: \Im\HH\rightarrow \Im \HH$ be a real
linear
traceless and bijectiv map which is symmetric with
respect to the standard scalar product $\ip_{\Bbb H}|_{\Im \Bbb H}$.
Then we define an inner product on $\Im \HH\oplus \Im \HH$ by
$$\langle (P_1,P_2), (Q_1,Q_2) \rangle := \langle P_1,AQ_2\rangle_{\Bbb H}+
\langle Q_1, A P_2\rangle_{\Bbb H} $$
and denote the resulting vector space with orthogonal quaternionic grading by
$\fa_A$.
We define $\alpha_+\in C^{2}(\HH^2,\fa_A)^{\Sp}$ by
$$ \alpha_+((p,q),(r,s)):=(\Im (\bar p s+\bar q r), \Im (\bar q s))\ .   $$
We claim that
\begin{equation}\label{hope}
    (\alpha_{+},0)\in  \cZ^{2}_{Q}(\HH^2,\Phi_{{\Bbb H}^2},\fa_{A})\ .
\end{equation}
We have to check that $\langle\alpha_+ \wedge \alpha_+\rangle =0 $. We
decompose $\bigwedge^{4}\HH^{2}={\displaystyle 
\bigoplus_{k+l=4}}\bigwedge^{k,l}$, where 
$\bigwedge^{k,l}=\bigwedge^{k}\HH\otimes\bigwedge^{l}\HH$. First we observe that
$\langle\alpha_{+} \wedge \alpha_{+}\rangle$ vanishes on
$\bigwedge^{k,l}$ for $(k,l)\ne (1,3)$. Next we
compute
\begin{eqnarray*}
 \textstyle{\frac{1}{2}}\langle\alpha_+ \wedge 
\alpha_+\rangle((1,0),(0,i),(0,j),(0,k))&=&
 \langle i, A(-i)\rangle_{\Bbb H}+\langle k, A(-k)\rangle_{\Bbb H}+\langle -j,A(j)\rangle_{\Bbb H}\\
 &=&-\tr A\,=\,0\ ,
 \end{eqnarray*}
and for $p,q\in \Im\HH\subset\HH$
 $$\textstyle{\frac{1}{2}}\langle\alpha_{+} \wedge
 \alpha_{+}\rangle((1,0),(0,p),(0,q),(0,1))\, = \,
 \langle p, A(-q)\rangle_{\Bbb H}+\langle -q, A(-p)\rangle_{\Bbb H}\,=\,0\ .
$$
$\Sp$-invariance implies that $\langle\alpha_{+} \wedge \alpha_{+}\rangle$ 
vanishes on
$\bigwedge^{1,3}$ as well. This proves (\ref{hope}).

Let $\fl_0,\fa_{0},\alpha_0,\gamma_0$ be as in Example~1. 
Fix $n\in \NN$ and choose a vector $\frak A=(A_1,\dots,A_n)$ of
traceless symmetric bijective maps $A_k:\Im\HH\rightarrow \Im \HH$. We equip
$$\fl:=\fl_{0}\oplus\HH^{n}\ ,\quad \fa_{\frak A}:=\fa_{0}\oplus
\bigoplus_{k=1}^{n} \fa_{A_k}$$
with their natural quaternionic gradings. Then we have natural 
projections
$$ \ph_{0}: \fl\longrightarrow \fl_{0}$$
and for $k=1,\dots,n$
$$ \ph_{k}: \fl\longrightarrow \HH^{2} ,\quad
\ph_{k}(p,P,p_{1},\dots,p_{n}):=(p,p_{k})\ .$$
We define
$$\alpha:=\ph_{0}^{*}\alpha_0\oplus\bigoplus_{k=1}^{n}\ph_{k}^{*}\alpha_+
\in C^{2}(\fl,\fa_{\frak A})^{\Sp} ,\quad \gamma:=\ph_{0}^{*}\gamma_0\in
C^{3}(\fl)\ .$$
Combining (\ref{hope}) with Lemma \ref{lemma} we see that
$(\alpha,\gamma)\in\cZ^{2}_{Q}(\fl,\Phi_{\fl},\fa_{\frak A})$.

\begin{lm}\label{emma}
 The cocycle $(\alpha,\gamma)\in \cZ^{2}_{Q}(\fl,\Phi_{\fl},\fa_{\frak A})$
 is admissible and indecomposable. 
\end{lm}

\proof
A straightforward verification shows that $(\alpha,\gamma)$ satisfies
the admissibility conditions of Definition~\ref{adm}.
In particular, $\alpha$ enjoys the following
properties
\begin{enumerate}
    \item[(a)] $\alpha(\bigwedge^{2}\fl)=\fa$,
    \item[(b)] For all $L\in\fl\setminus\{0\}$ there exists $L'\in\fl$ 
    such that $\alpha(L,L')\ne 0$.
\end{enumerate}
Now we use Properties (a) and (b) in order to show indecomposability 
of $(\alpha,\gamma)$. Assume that we have decompositions
$$ \fl=\fl_{1}\oplus \fl_{2}, \ \fa=\fa_{1}\oplus \fa_{2},\ 
\alpha=\pro_{1}^{*}\alpha_{1}\oplus \pro_{2}^{*}\alpha_{2}\ , $$
for certain $\alpha_{i}\in C^{2}(\fl_{i},\fa_{i})^{\Sp}$. We may 
assume that $\fl_{1}$ is non-abelian. Then $\Sp$-invariance of 
$\fl_{1}$ implies that $\fl_{+}=(\fl_{0})_{+}=(\fl_{1})_+$.
This forces $(\fl_{2})_+=[(\fl_2)_-,(\fl_2)_-]=0$, which in turn implies
$\fl_{2}\subset \HH^{n}$. It follows that 
$\fa_{2}=\alpha(\bigwedge^{2}\fl_{2})$ is isotropic, hence 
$\fa_{2}=\{0\}$, $\alpha_{2}=0$. Therefore $\alpha=\pro_{1}^{*}\alpha_{1}$.
Now (b) implies that $\fl_{2}=\{0\}$. It follows that $(\alpha,\gamma)$ 
is indecomposable.
\qed

As in Example 1 we obtain
\begin{co} 
$\fd_{\alpha,\gamma}(\fl,\Phi_{\fl},\fa_{\frak A})$ is an indecomposable
hyper-K\"ahler symmetric triple of signature $(4n+4,4n+12)$. The holonomy group
of a  symmetric space associated with this triple is non-abelian.
\end{co}

\subsection*{Example 3}

Let $(\fl,\Phl)$ be a Lie algebra with proper quaternionic grading.
Then we can form the cotangent Lie algebra $T^*\fl:=\fl\ltimes\fl^*$
which possesses a quaternionic grading $T^*\Phi_\fl$ and an invariant metric $\ip_{T^*\fl}$ (given by the dual pairing) in a natural way. This construction could be viewed
as a special case of the one presented in Section \ref{SCon}
-- apart from the fact that the involved cocycle $(0,0)$ is not admissible.
If we require in addition that
\begin{equation}\label{warenhaus}
\fz(\fl)\subset \fl_-\ ,
\end{equation}
then $(T^*\fl,T^*\Phl,\ip_{T^*\fl})$ is a
hyper-K\"ahler symmetric triple.

If $(\fg,\Phg,\ip_\fg)$ is a hyper-K\"ahler symmetric triple, then $\fg$
satisfies (\ref{warenhaus}), and the cotangent hyper-K\"ahler symmetric triple $(T^*\fg,T^*\Phg,\ip_{T^*\fg})$ is isomorphic to the tangent triple
$(T\fg,T\Phg,\ip_{T\fg})$, where
$$ T\fg:=\fg\ltimes\fg\, ,\quad \langle (X_1,Y_1), (X_2,Y_2)\rangle_{T\fg}
:=\langle X_1,Y_2\rangle_{\fg}+\langle X_2,Y_1\rangle_{\fg}\ .$$

\begin{pr}\label{uff}
\begin{enumerate}
\item[(a)] Let $(\fl,\Phl)$ be a non-abelian Lie algebra with proper quaternionic grading
satisfying (\ref{warenhaus}). If $(\fl,\Phl)$ is indecomposable as a Lie algebra
with quaternionic grading, then the hyper-K\"ahler symmetric triple
$(T^*\fl,T^*\Phl,\ip_{T^*\fl})$ is indecomposable.
\item[(b)] Let $(\fg,\Phg,\ip_\fg)$ be an indecomposable hyper-K\"ahler
symmetric triple with $\fg$ non-abelian. Then the hyper-K\"ahler symmetric triple
$(T\fg,T\Phg,\ip_{T\fg})$ is indecomposable as well.
\end{enumerate}
\end{pr}
\proof
Let $(\fg,\Phg,\ip_\fg)$ be an indecomposable hyper-K\"ahler
symmetric triple with $\fg$ non-abelian. By Lemma \ref{A0} the tuple
$(\fg,\ip_\fg)$ is indecomposable as a metric Lie algebra. Then, according to \cite{A78}, Theorem 5, the metric Lie algebra $(T\fg,\ip_{T\fg})$ is indecomposable as well.
This implies $(b)$.

Now let $(\fl,\Phl)$ be as in $(a)$. We consider the Lie algebra $\fsp(1)\ltimes\fl$,
where $\fsp(1)$ acts on $\fl$ by the differential of $\Phl$. We let $\fsp(1)\ltimes\fl$ act on $T^*\fl$
by $(Q,L)v:=d(T^*\Phl)(Q)v+[L,v]$, $Q\in\fsp(1)$, $L\in\fl$, $v\in T^*\fl$.
Then $T^*\fl=\fl\oplus\fl^*$ is a decomposition into indecomposable $\fsp(1)\ltimes\fl$-submodules. Let us assume that there is a non-trivial
decomposition into hyper-K\"ahler symmetric
triples $T^*\fl=\fh_1\oplus \fh_2$. In particular, $\fh_1$ and $\fh_2$ are
$\fsp(1)\ltimes\fl$-submodules. The Krull-Schmidt Theorem (see e.g.~\cite{Jac}, p.~115) implies
that one of these modules, say $\fh_1$, is isomorphic to $\fl$. Pulling
back the invariant bilinear form $\ip_{T^*\fl}|_{\fh_1}$ from $\fh_1$ to $\fl$
we obtain an inner product $\ip_\fl$ such that $(\fl,\Phl,\ip_\fl)$ is
an indecomposable hyper-K\"ahler symmetric triple. Now we apply $(b)$
to this triple which yields a contradiction to the assumption that
$T^*\fl\cong T\fl$ is decomposable as a hyper-K\"ahler symmetric triple.
This finishes the proof of $(a)$.
\qed

If the symbol $\fd$ stands for a hyper-K\"ahler symmetric triple, then we
simply write $T\fd$ or $T^1\fd$ for its tangential triple. We form
higher tangential triples by $T^n\fd:=T(T^{n-1}\fd)$, $n\ge 2$.

\begin{co}\label{Tang}
Let $\fd$ be one of the hyper-K\"ahler symmetric triples constructed
in Examples 1 and 2 of signature $(4k,4k+8)$, $k\ge 1$. Then $T^n\fd$,
$n\ge 1$, is an indecomposable hyper-K\"ahler symmetric triple of
neutral signature $(2^{n+2}(k+1),2^{n+2}(k+1))$. 
The holonomy group
of a  symmetric space associated with $T^n\fd$ is non-abelian.
\end{co}

\begin{re}
{\rm
If $\fd=\dd$ for some admissible cocycle $(\alpha,\gamma)\in \cZP$ (as e.g. in Corollary \ref{Tang}), then
$T\fd\cong\fd_{T\alpha,T\gamma}(T\fl,T\Phi_{\fl},T\fa)$, where
\begin{eqnarray*}
T\alpha\left((L_1,L'_1),(L_2,L'_2)\right)&=&
\left(\alpha(L_1,L_2),\alpha(L_1,L'_2)+\alpha(L'_1,L_2)\right)\ ,\\
T\gamma\left((L_1,L'_1),(L_2,L'_2),(L_3,L'_3)\right)&=&
\gamma(L_1,L_2,L_3')+\gamma(L_1,L'_2,L_3)+\gamma(L'_1,L_2,L_3)\ .
\end{eqnarray*}
Moreover, one can show that the cocycle $(T\alpha,T\gamma)$ is admissible
again. This
last assertion is not true for the general symmetric triples
$\fd_{\alpha,\gamma}(\fl,\theta_{\fl},\fa)$ constructed in \cite{KO3}
from admissible
cocycles $(\alpha,\gamma)\in\cZ^{2}_{Q}(\fl,\theta_{\fl},\fa)$, only for those with nilpotent $\fl$ and a trivial $\fl$-module $\fa$.}
\end{re}

\section{A classification scheme} 
\label{SClass} 
Let $(\fl,\Phl)$ be a Lie algebra with proper quaternionic grading and let
$(\fa,\ip_\fa,\Pha)$ be a vector space with orthogonal quaternionic grading.
We consider their automorphism groups $\Aut(\fl,\Phl)$, and 
$\Aut(\fa,\ip_{\fa},\Pha)$ and form 
$$G_{\fl,\Phl,\fa}:=\Aut(\fl,\Phl)\times\Aut(\fa,\ip_{\fa},\Pha)\ .$$ 
There is a natural right action of $G_{\fl,\Phl,\fa}$ on $\cHP$ which leaves 
$\cHP_0$ invariant. 
\begin{theo} \label{general}
 There is a bijective map from the set of isomorphism classes of 
 indecomposable hyper-K\"ahler symmetric triples 
 to the union of orbit spaces 
 $$\coprod_{(\fl,\Phi_{\fl})} \coprod_{(\fa,\ip_{\fa},\Pha)} 
 \cHP_{0}/ G_{\fl,\Phl,\fa},$$ 
where the union is taken over all isomorphism classes of Lie algebras 
$(\fl,\Phi_{\fl})$ with proper quaternionic 
grading and all isomorphism classes of vector spaces $(\fa,\ip_{\fa},\Pha)$ 
with orthogonal quaternionic grading. 

The inverse of this map sends the orbit of $[\alpha,\gamma]\in\cHP_{0}$ 
to the isomorphism class of $\dd$. 
\end{theo} 

\proof 
The theorem is the hyper-K\"ahler analog of our classification scheme of 
symmetric triples \cite{KO3}, Theorem 6.1. The proof of the latter is 
contained in \cite{KO3}, Sections 4~--~6. It carries over to the present 
situation, 
one has to take care of the quaternionic gradings, only. 
For the convenience of the reader we describe the map which associates 
a cohomology class to a given hyper-K\"ahler symmetric triple 
$(\fg,\Phi_\fg,\ip_\fg)$.

By Corollary \ref{nico} the Lie algebra $\fg$ is nilpotent.
We look at its lower central series
and form the isotropic ideal 
$$\fri:= \sum_{k=2}^{\infty} \fg^k\cap(\fg^k)^\perp\ .$$ 
Then $\fl:=\fg/\fri^\perp$ and $\fa:=\fri^\perp/\fri$ inherit the desired 
structures from $\fg$. Moreover, the nilpotency of $\fg$ implies that the induced
$\fl$-action on $\fa$ is trivial. We choose an $Sp(1)$-equivariant section
$s:\fl\rightarrow \fg$ with isotropic image and set
\begin{eqnarray*} 
\alpha(L_1,L_2)&:=&[s(L_1),s(L_2)]_\fg-s([L_1,L_2]_\fl)\ \mod \fri\quad 
\\ 
\gamma(L_1,L_2,L_3)&:=& \langle\, [s(L_1),s(L_2)]_\fg, s(L_3)\rangle_\fg \ . 
\end{eqnarray*} 
Then $(\alpha,\gamma)$ is an admissible cocycle in $\cZP$.
It is indecomposable if $(\fg,\Phi_\fg,\ip_\fg)$ is so. The desired cohomology class is given by
$[\alpha,\gamma]\in\cHP_0$.
\qed 

\section{Classification in the case of signature $(4,4n)$} \label{S4}

Our general classification scheme can be used to find explicit classification results (i.e.~lists), if we only consider hyper-K\"ahler symmetric spaces of a
given small index. Here we will classify indecomposable
hyper-K\"ahler symmetric spaces of index~4. 
We want to apply Theorem \ref{general}. So we have to decide for which 
$(\fl,\Phl)$ 
there is a vector space with orthogonal quaternionic grading $\fa$ and an 
element $[\alpha,\gamma]\in  \cHP_{0}$ such that $\dd$ has index~$4$.

\begin{pr} \label{P1} Let $(\fl,\Phl)$ be a Lie algebra with proper quaternionic 
grading, let 
$\fa$ be a vector space with orthogonal quaternionic grading and suppose 
$[\alpha,\gamma]\in\cHP_0$. If $\dd$ is a hyper-K\"ahler symmetric triple of 
index 4, then $\fl=0$ or $(\fl,\Phl)$ is isomorphic either to 
$(\fl_0,\Phi_{\fl_0})$ (as defined in Section \ref{SEx}) or to $(\HH,\Phi_{\Bbb 
H})$, where 
$\Phi_{\Bbb H}$ is the left multiplication on~$\HH$. 
\end{pr} 
\proof 
{}From (\ref{signature}) we know that  $\dim\fl_{-}\in\{0,4\}$. If 
$\dim \fl_{-}=0$, then $\fl=0$. Suppose now that $\dim \fl_{-}=4$,  
i.e.~$\fl_{-}\cong\HH$. Consider now the map $\ad L|_{\fl_{-}}$ for an arbitrary  
$L\in\fl_+$. This map 
commutes with the $\Sp$-action, i.e.~it is $\HH$-linear. On the other 
hand it must be nilpotent, since $\fl$ is nilpotent by Proposition 
\ref{nilpferd}. Hence the map is zero, and it follows that 
$[\fl_{-},\fl_{+}]=0$. This gives 
$[\fl_{+},\fl_{+}]=[[\fl_{-},\fl_{-}],\fl_{+}]=0$. In particular, we 
have $\fl'=\fl_{+}\subset\fz(\fl)$. 

Since 
$\Sp$ acts by automorphisms we have 
$$[1,i]=-[j,k],\quad [1,j]=[i,k],\quad [1,k]=-[i,j].$$
In particular we get $\dim \fl'\le 3$. If $\dim \fl'=0$, then 
$\fl\cong\HH$, if $\dim \fl'=3$, then $\fl\cong\fl_{0}$. Assume $\dim 
\fl'=1$. By an easy direct computation or using \cite{KO2}, Prop.~6.2.~one 
obtains
$\cHP_{\sharp}=\emptyset$, a contradiction. 

It remains to exclude the case $\dim \fl'=2$. If $\dim \fl'=2$, then we may 
assume 
$$ [1,i]=-[j,k]=:X,\quad [1,j]=[i,k]=:Y,\quad [1,k]=-[i,j]=0.$$
Then $d\alpha=0$ implies $\alpha(X,Y)=0$ and 
$0=d\alpha(1,i,j)=\alpha([1,i],j)+\alpha([j,1],i)$, thus 
$\alpha(X,j)=\alpha(Y,i)=:A.$
This together with the $\Sp$-invariance of $\alpha$ implies 
$\fa_{-}=\alpha(\fl_{-},\fl_{+})=\Span\{A,iA,jA,kA\}$. But 
$d\gamma(i,j,X,Y)=0$ now implies 
$$0=\langle \alpha(j,X),\alpha(i,Y)\rangle +\langle 
\alpha(X,i),\alpha(j,Y)\rangle=2\langle A,A \rangle .$$ 
Hence $A$ is isotropic. Since $\fa_{-}$ is non-degenerate it follows 
that $A=0$, thus $\alpha(\fl_{-},\fl_{+})=0$. In particular, 
$\alpha(X,\fl)=0$. Hence $\alpha$ satisfies Condition $(i)$ of $(A_{1})$ 
for $L_{0}=X$. Let us show that Condition $(ii)$ is also satisfied. 
Because of the $\Sp$-invariance of $\gamma$ we have 
$\gamma(\fl_{+},\fl_{+},\fl_{-})=0$. Since, moreover, $\dim 
\fl_{+}=2$ we obtain $\gamma(\fl_{+},\fl_{+},\fl)=0$. In particular,  
$\gamma(L,X,\cdot)=0\in(\fl')^{*}$ for all $L\in\fl$. Hence $(ii)$ is 
satisfied for $A_{0}=0$ and $Z_{0}=0$. Now $(A_{1})$ implies $X=0$, a 
contradiction. 
\qed 

Next we have to determine $\cHP_{0}/G_{\fl,\Phl,\fa}$ for all $(\fl,\Phl)$ 
which appear in Proposition \ref{P1} and for all (suitable) $\fa$. 
The case $\fl=0$ is trivial. Thus let
us start with 
$(\fl,\Phl)=(\HH,\Phi_{\Bbb H})$. 
For a fixed orthonormal basis $A_1,A_2$ of $\fa=\fa_+=\RR^{1,1}$ we define 
$\alpha'\in Z^2(\HH,\fa)^{\Sp}$ by  
$$\alpha'(1,i)=A_1,\quad \alpha'(1,j)=A_2,\quad \alpha'(1,k)=0.$$ 

Furthermore, for a fixed orthonormal basis $A_1,A_2,A_3$ of 
$\fa=\fa_+=\RR^{1,2}$ and each real number $0<r\le\pi/4$ we define $\alpha_r\in 
Z^2(\HH,\fa)^{\Sp}$ by 
$$\alpha_r(1,i)=A_1,\quad \alpha_r(1,j)=\sin r\cdot A_2,\quad 
\alpha_r(1,k)=\cos r \cdot A_3.$$ 

Analogously, for a fixed orthonormal basis $A_1,A_2,A_3$ of 
$\fa=\fa_+=\RR^{2,1}$ and each real number $0<s\le\pi/4$ we define $\alpha_s\in 
Z^2(\HH,\fa)^{\Sp}$ by 
$$\alpha_s(1,i)=\sin s\cdot A_1,\quad \alpha_s(1,j)=\cos s \cdot A_2,\quad 
\alpha_s(1,k)=A_3.$$ 

Then it is easy to verify that $\langle \alpha'\wedge\alpha'\rangle=0$, 
$\langle \alpha_r\wedge\alpha_r\rangle=0$, and $\langle 
\alpha_s\wedge\alpha_s\rangle=0$ holds, e.g. 
\begin{eqnarray*} 
\langle \alpha_r\wedge\alpha_r\rangle(1,i,j,k)&=&\langle 
\alpha_r(1,i),\alpha_r(j,k)\rangle+\langle \alpha_r(j,1),\alpha_r(i,k)\rangle+ 
\langle\alpha_r(i,j) ,\alpha_r(1,k)\rangle\\ 
&=&\langle A_1,-A_1\rangle+\langle -\sin r\cdot A_2,\sin r\cdot A_2\rangle+ 
\langle -\cos r \cdot A_3 ,\cos r \cdot A_3\rangle\\&=&1-\sin^2r-\cos^2r\ =\ 0. 
\end{eqnarray*} 
\begin{pr} \label{P2} 
We equip the spaces $\RR^{1,1}$, $\RR^{1,2}$ and $\RR^{2,1}$ with the trivial $\Sp$-action. Then we have 
\begin{eqnarray*} 
\cH^2_Q(\HH,\Phi_{\Bbb H},\RR^{1,1})_0/G_{{\Bbb H},\Phi_{\Bbb H},{\Bbb R}^{1,1}} 
&=& \{[\alpha',0]\}\\ 
\cH^2_Q(\HH,\Phi_{\Bbb H},\RR^{1,2})_0/G_{{\Bbb H},\Phi_{\Bbb H},{\Bbb R}^{1,2}} 
&=& \{[\alpha_r,0] \mid 0<r\le\pi/4 \}\\ 
\cH^2_Q(\HH,\Phi_{\Bbb H},\RR^{2,1})_0/G_{{\Bbb H},\Phi_{\Bbb H},{\Bbb R}^{2,1}} 
&=& \{[\alpha_s,0] \mid 0<s\le\pi/4 \}. 
\end{eqnarray*} 
If $\fa$ is a vector space with orthogonal quaternionic grading which is not 
isomorphic to $\RR^{1,1}$, $\RR^{1,2}$ or $\RR^{2,1}$ with the trivial $\Sp$-action, then 
$\cH^2_Q(\HH,\Phi_{\Bbb H},\fa)_0$ is empty. 
\end{pr}
\proof Let $\fa$ be a vector space with orthogonal quaternionic grading. 
Since $C^3(\HH)^{\Sp}$ vanishes and $H^2(\HH,\fa)=C^2(\HH,\fa)$ we have
$$\cH^2_Q(\HH,\Phi_{\Bbb H},\fa)=\{[\alpha,0]\mid \alpha\in C^2(\HH,\fa)^{\Sp},\ 
\langle \alpha\wedge \alpha\rangle=0\}.$$ 
Each element $\alpha\in C^2(\HH,\fa)^{\Sp}$ is given by the values 
\begin{equation}\label{AAA} 
\alpha(1,i)=A'_1,\ \alpha(1,j)=A'_2,\ \alpha(1,k)=A'_3,\quad A'_i\in\fa_+,
\end{equation} 
and $\langle \alpha\wedge \alpha\rangle=0$ holds if and only if $\langle 
A_1',A'_1\rangle+ \langle A'_2,A'_2\rangle+ \langle A'_3,A'_3\rangle=0.$ 
Moreover, 
$[\alpha,0]$ is indecomposable if and only if $\fa=\fa_+=\Span\{A'_1,A'_2,A'_3\}$. For admissibility we have to check only $(A_0)$ and $(B_0)$. For an indecomposable cohomology class $[\alpha,0]$ these conditions are satisfied if and only if $\alpha\not=0$.
In particular, if $\cH^2_Q(\HH,\Phi_{\Bbb H},\fa)_0$ is not empty, then $\fa=\fa_+$ is 
isomorphic to $\RR^{1,1}$, $\RR^{1,2}$ or $\RR^{2,1}$. Now we study the action of 
$G_{{\Bbb H},\Phi_{\Bbb H},\fa}$ on $\cH^2_Q(\HH,\Phi_{\Bbb H},\fa)_0$. We have 
$\Aut(\HH,\Phi_{\Bbb H})\cong \HH^*$, where an element $q\in \HH^*$ acts by right 
multiplication (denoted by $R_q$) on $\HH$. 
Take $q\in \HH^*$. Then $q=rq_0$ with $r\in\RR,r>0$ and $q_0\in\Sp$. Take 
$\alpha$ as in (\ref{AAA}) and set 
$$R_q^*\alpha(1,i)=\bar A_1,\ R_q^*\alpha(1,j)=\bar A_2,\ R_q^*\alpha(1,k)=\bar 
A_3.$$ 
Let $\lambda:\Sp\rightarrow \SO(3)$ be the double covering. Then we have 
$$(\bar A_1,\bar A_2,\bar A_3)=r^2(A'_1,A'_2,A'_3)\cdot \lambda(\bar q_0).$$ Hence, 
replacing $\alpha$ by an element in the same $G_{{\Bbb H},\Phi_{\Bbb 
H},\fa}$-orbit we may assume that $A'_1,A'_2,A'_3$ are pairwise orthogonal and 
that 
$$\langle A'_1,A'_1 \rangle =-1,\ \langle A'_2,A'_2 \rangle=1,\ A'_3=0$$ 
or 
$$\langle A'_1,A'_1 \rangle =-1,\ \langle A'_2,A'_2 \rangle=t, \langle A'_3,A'_3 
\rangle=1-t,\ 0<t\le 1/2 $$ 
or 
$$\langle A'_1,A'_1 \rangle =-t,\ \langle A'_2,A'_2 \rangle=t-1, \langle 
A'_3,A'_3 \rangle=1,\ 0<t\le 1/2, $$ 
which immediately implies the assertion. 
\qed

Now let us consider $(\fl_0,\Phi_{\fl_{0}})$. Let $\fa_{0}$  and
$[\alpha_0,\gamma_0]\in \cH^2_Q(\fl_0,\Phi_{\fl_0},\fa_{0})_0$ be 
as defined in Section \ref{SEx}. 
 
\begin{pr} \label{P3}
The orbit space 
$\cH^2_Q(\fl_0,\Phi_{\fl_0},\fa_{0})_0/G_{\fl_0,\Phi_{\fl_0},\Phi_{\fa_{0}}}$ 
contains 
exactly one element (represented by $[\alpha_0,\gamma_0]$). If $\fa$ is a 
vector space with orthogonal quaternionic grading which is not isomorphic to 
$\fa_{0}$ or to $\fa_{0}^{-}:=(\fa_{0},-\ip_{0},\Phi_{\fa_{0}})$, then 
$\cH^2_Q(\fl_0,\Phi_{\fl_0},\fa)_0$  is empty. 
\end{pr}
\proof Let $\fa$ be a vector space with orthogonal quaternionic grading. It is easy to verify that for 
$$ Z_{\fl_0}^2:=\left\{\alpha\in C^2(\fl_0,\fa)^{\Sp} \ \Big| 
\begin{array}{l} i\alpha(1,I)+j\alpha(1,J)+k\alpha(1,K)=0\\ 
\alpha((\fl_0)_+,(\fl_0)_+)=\alpha((\fl_0)_-,(\fl_0)_-)=0 
\end{array}\right\}$$ 
the map 
$$ Z_{\fl_0}^2 \longrightarrow H^2(\fl_0,\fa)^{\Sp},\quad \alpha\longmapsto 
[\alpha]$$ 
is correctly defined (similar computation as in the proof of Lemma \ref{lemma}) and an isomorphism.  
Now assume that $\cH^2_Q(\fl_0,\Phi_{\fl_0},\fa)_0$ is not empty and take 
$[\alpha,\gamma]\in\cH^2_Q(\fl_0,\Phi_{\fl_0},\fa)_0$. We may assume that 
$\alpha\in Z_{\fl_0}^2$. We set 
$$A_1:=-\alpha(i,I),\quad A_2:=-\alpha(j,J),\quad A_3=-\alpha(k,K)$$ and 
$\gamma(I,J,K)=2c$. 
Since $(\alpha,\gamma)$ is indecomposable we have $$\fa=\fa_-=\HH \otimes_{\Bbb 
R}\Span\{A_1,A_2,A_3\}.$$ By the same computation as in the proof of Lemma 
\ref{lemma} the equation 
$$\textstyle{\frac 12} \langle \alpha\wedge 
\alpha\rangle(1,q,P,Q)=d\gamma(1,q,P,Q)$$ 
for all $q\in \HH$, $P,Q\in\{I,J,K\}$ yields 
$$\langle q_1 A_1, q_2A_2\rangle_\fa = \langle q_1 A_1, q_2A_3\rangle_\fa 
= \langle q_1 A_2, q_2A_3\rangle_\fa = -\textstyle{\frac 12}\langle q_1, 
q_2\rangle_{\Bbb H}\cdot \gamma(I,J,K)=- c\langle q_1, q_2\rangle_{\Bbb H}.$$ 
This implies 
$$0=\langle p(A_1+A_2+A_3), qA_i\rangle=\langle pA_i, qA_i\rangle_\fa-2c\langle 
p, q\rangle_{\Bbb H},$$ 
thus $\langle pA_i, qA_i\rangle_\fa=2c\langle p, q\rangle_{\Bbb H}$ for 
$i=1,2,3$, $p,q\in\HH$. Assume $c=0$. 
Then $\alpha=0$ and $\gamma((\fl_{0})_{+},(\fl_{0})_{+},\fl_{0})=0$. 
Then it follows in the same way as in the proof of Proposition~\ref{P1}, that $[\alpha,\gamma]$ is not admissible. Hence $c\not=0$.   If $c>0$, 
then by the above equations $\fa$ is 
isomorphic to $\fa_0$ as a vector space with orthogonal quaternionic grading and 
if $c<0$, then $\fa$ is isomorphic to $\fa_0^{-}$. 

Hence we may assume that 
$\fa=\fa_{0}$, $\alpha=|c|^{1/2}\alpha_0$, 
or $\fa=\fa_{0}^{-}$, $\alpha=|c|^{1/2}\alpha_0$. Let us now show 
that $[\alpha,\gamma]=[\alpha,\tilde\gamma]\in\cHP_{0}$ with 
$\tilde \gamma(I,J,K)=2c$ and $\tilde 
\gamma((\fl_{0})_{-},(\fl_{0}),(\fl_{0}))=0$. We consider the subspaces
$$\{\langle \alpha\wedge\tau\rangle \in C^3(\fl_0)^{\Sp}\mid \tau \in 
Z^{1}(\fl_{0},\fa)^{\Sp}\}\subset \{\gamma'\in 
Z^3(\fl_0)^{\Sp} \mid \gamma'(I,J,K)=0\}$$ of $C^3(\fl_0)^{\Sp}$.
By definition of $\cHP$ it suffices to show that both subspaces are 
equal. Note that the latter of these spaces is $8$-dimensional. Since 
$Z^1(\fl_{0},\fa)^{\Sp}$ is also $8$-dimensional it remains to prove 
that the map 
$$Z^{1}(\fl_{0},\fa)^{\Sp}\ni \tau \longmapsto \langle 
\alpha\wedge\tau\rangle \in C^3(\fl_0)^{\Sp}$$ 
is injective. Assume that $\tau$ is in the kernel of this map. Then we 
have $$0=\langle \alpha\wedge\tau\rangle(1,q,Q)=\langle  
\alpha(q,Q),\tau(1)\rangle-\langle  
\alpha(1,Q),\tau(q)\rangle= 2\langle  
\alpha(q,Q),\tau(1)\rangle$$
for all $q\in\Im \HH$ and all $Q\in\{I,J,K\}$. Since 
$$\fa=\Span\{\alpha(q,Q)\mid q\in\Im \HH, Q\in\{I,J,K\}\}$$ this implies 
$\tau(1)=0$, hence $\tau=0$.

Hence we may assume that 
$\fa=\fa_{0}$, $\alpha=|c|^{1/2}\alpha_0$, and $\gamma= |c|\gamma_0$ 
or $\fa=\fa_{0}^{-}$, $\alpha=|c|^{1/2}\alpha_0$, and $\gamma= -|c|\gamma_0$. If 
we 
now take $S=|c|^{-1/3}\Id \oplus |c|^{-1/6}\Id:(\fl_0)_+\oplus 
(\fl_0)_-\rightarrow (\fl_0)_+\oplus (\fl_0)_-$, then $S\in 
\Aut(\fl_0,\Phi_{\fl_0})$ and we get $(S^*\alpha, 
S^*\gamma)=(\alpha_0,\pm\gamma_0)$. 
\qed 

As a consequence of Theorem \ref{general}, Equation (\ref{signature}) 
and Propositions \ref{P1} -- \ref{P3} we obtain the following 
classification.
\begin{theo} \label{T4}
If $(\fg,\Phi_\fg,\ip_\fg)$ is a hyper-K\"ahler symmetric triple which is 
associated with an indecomoposable hyper-K\"ahler symmetric space of index 4, 
then it is isomorphic to $(\HH,\Phi_{\Bbb H},-\ip_{\Bbb H})$ or to $\dd$ for 
exactly one of the data in the following table: 
\begin{center} 
\begin{tabular}{|c|c|c|c|c|} 
\hline 
&&&&\\[-2ex] 
$(\fl,\Phl)$ & $\fa$ & $\alpha$ &$\gamma$ & parameters \\[1ex] 
\hline \hline 
&&&&\\[-2.0ex] 
$(\HH,\Phi_{\Bbb H})$ & $\fa=\fa_+=\RR^{1,1}$ & $\alpha'$ & $0$ & -- 
\\[0.5ex]\cline{2-5}&&&& \\[-2ex] 
& $\fa=\fa_+=\RR^{1,2}$ & $\alpha_r$ & $0$ & $0<r\le\pi/4 $ 
\\[0.5ex]\cline{2-5}&&&& \\[-2ex] 
& $\fa=\fa_+=\RR^{2,1}$ & $\alpha_s$ & $0$ & $0<s\le\pi/4 $ 
\\[0.5ex]\hline&&&& \\[-2ex] 
$(\fl_0,\Phi_{\fl_0})$ & $\fa_0$ & $\alpha_0$ & $\gamma_0$ & -- 
\\[0.5ex]\hline 
\end{tabular} 
\end{center} 
\end{theo}

\vspace{1.6cm}
{\footnotesize 

Ines Kath\\ 
Max-Planck-Institut f\"ur Mathematik in den Naturwissenschaften\\ 
Inselstra{\ss}e 22-26, D-04103 Leipzig, Germany\\ 
email: ikath@mis.mpg.de\\[2ex] 
Martin Olbrich\\ 
Mathematisches Institut der Georg-August-Universit{\"a}t G{\"o}ttingen\\ 
Bunsenstr. 3-5, D-37073 G{\"o}ttingen, Germany\\ 
email: olbrich@uni-math.gwdg.de} 

\end{document}